\documentclass[journal]{IEEEtran}
\usepackage{amsmath, amssymb}
\newtheorem{theorem}{Theorem}[section]

\newtheorem{corollary}[theorem]{Corollary}
\newtheorem{lemma}[theorem]{Lemma}
\newtheorem{observation}[theorem]{Observation}

\newtheorem{definition}[theorem]{Definition}

\numberwithin{equation}{section}

\def \endproof {{\mbox{}\nolinebreak\hfill\rule{2mm}{2mm}\par\medbreak}}

\newlength{\algorithmwidth}
\def \R {\mathbb{R}}
\def \C {\mathbb{C}}

\def \E {\mathbb{E}}
\def \P {\mathbb{P}}

\def \one {{\bf 1}}

\def \NN {\mathcal{N}}

\def \OO {\mathcal{O}}

\def \a {\alpha}
\def \b {\beta}

\def \e {\varepsilon}

\def \d {\delta}

\def \l {\lambda}

\def \s {\sigma}
\def \t {\tau}

\def \< {\langle}
\def \> {\rangle}
\def \^ {\widehat}

\def \id {{\it id}}

\def \supp {{\rm supp}}

\begin{document}
\title {Uncertainty principles and vector quantization}

\author{Yurii~Lyubarskii
  and Roman~Vershynin%
  \thanks{Yu. Luibarskii in with Norwegian University of Science and Technology. 
  He was partially supported by the Norwegian Research Council projects 160192/v30, 10323200.}
  \thanks{R. Vershynin is with University of Michigan.
  He was partially supported by NSF grants DMS 0401032, 0652617, 0918623 and Alfred P. Sloan Foundation.}
}

\maketitle

\begin{abstract}
 Given a frame in $\C^n$ which satisfies a form of the uncertainty principle
 (as introduced by Candes and Tao), it is shown how to 
 quickly convert the frame representation of every vector
 into a more robust Kashin's representation whose coefficients
 all have the smallest possible dynamic range $O(1/\sqrt{n})$.
 The information tends to spread evenly among these coefficients.
 As a consequence, Kashin's representations have a great power for reduction
 of errors in their coefficients, including coefficient losses and distortions.
\end{abstract}

\begin{keywords}
Frame representations, Kashin's representations, restricted isometries, uncertainty principles
\end{keywords}

\section{Introduction}
\label{s: acg}

Quantization is a representation of continuous structures with discrete structures.
Digital signal processing, which has revolutionized the modern
treatment of still images, video and audio, employs quantization as
a conversion step from the analog to  digital world.
A survey of the state-of-the-art of quantization prior to 1998
as well as outline of its numerous applications can be found the
 paper \cite{GN} by Gray and Neuhoff. For more recent developments,
we refer the reader to \cite{DY} and  references therein.

In this paper, we are interested in robust vector encoding and vector quantization.
Orthogonal expansions gives a classical way to encode vectors in finite dimensions.
One first chooses a convenient orthonormal basis $(u_i)_{i=1}^n$ of $\C^n$.
Then one encodes a vector $x \in \C^n$ by the coefficients $(a_i)_{i=1}^n$
of its orthogonal expansion
$$
x = \sum_{i=1}^n a_i u_i, \quad \text{where } a_i = \< x, u_i \> .
$$
An example of this situation is the discrete Fourier transform.
At the next step, one quantizes the coefficients $a_i$ using a convenient
scalar quantizer (for example, a uniform quantizer with fixed number of levels).

A drawback of orthogonal expansions is that the information contained
in the vector $x$ may get distributed unevenly among the coefficients $a_i$,
which makes this encoding vulnerable to distortions and losses of the coefficients.
For example, if $x$ is collinear with the first basis vector $u_1$ then all the
coefficients except $a_1$ are zero. If the first coefficient $a_1$ is lost
(for example due to transmission failure) then we can not reconstruct
the vector $x$ even approximately.

A popular way to improve the stability of vector encoding
is to use redundant systems of vectors $(u_i)_{i=1}^N$ in $\C^n$ called tight frames.
These are generalizations of orthonormal bases in the sense that every vector
$x \in \C^n$ can still be represented as
\begin{equation}                            \label{frame}
 x = \sum_{i=1}^N a_i u_i, \quad \text{where } a_i = \< x, u_i \> ,
\end{equation}
but for $N > n$ frames are clearly linearly dependent systems of vectors.
These dependencies cause the information contained in $x$ to spread
among several frame coefficients $a_i$, which improves the stability
of such representations with respect to errors (for example losses and quantization errors),
see e.g. \cite{D, GVT, C} and references therein.

The idea of spreading the information evenly among the coefficients is developed
in the present paper, and in a sense it is pushed to its limit.
As in the previous approaches, we shall start with a frame $(u_i)_{i=1}^N$.
But instead of the standard frame expansions \eqref{frame}
we will be looking at expansions $x = \sum_{i=1}^N a_i u_i$ with coefficients having the smallest
possible dynamic range $|a_i| = O(1/\sqrt{N})$.
This ensures that the information contained in $x$ is spread among the coefficients $a_i$
nearly uniformly. We call such representations of vector $x$ {\em Kashin representations}.
In this paper we demonstrate the following:

(a) there exist frames $(u_i)_{i=1}^N$ in $\C^n$ with redundancy factor $N/n$
as close as one likes to $1$, and such that every vector $x \in \C^n$ has a
Kashin representation;

(b) such frames are those that satisfy a form of the uncertainty principle.
More precisely, their matrices satisfy a weak version of the restricted isometry
property introduced by Candes and Tao \cite{CT2}.
In particular, many natural random frames have this property;

(c) there is a fast algorithm which converts  frame representation \eqref{frame}
into a Kashin representation of $x$.

Kashin's representations withstand errors in their coefficients in a very strong way --
the  representation error gets bounded by the {\em average}, rather than the sum,
of the errors in the coefficients.
These errors may be of arbitrary nature, including distortion (e.g. due to scalar quantization)
and losses (e.g. due to transmission failures).

The article is organized as follows.
Section~\ref{s: Kashin} introduces Kashin's representations,
discusses their relation to convex geometry (Euclidean projections of the cube)
and explains how one can use Kashin's representations for vector quantization.
In Section~\ref{s: computing Kashin}, we discuss the uncertainty principle
for matrices and frames. Theorems~\ref{algorithm valid} and \ref{almost tight} state
that for frames that satisfy the uncertainty principle, every frame representation
can be replaced by  Kashin's representation. A robust algorithm is given
to quickly convert frame into Kashin's representations.
In Section \ref{s: matrices UP}, we discuss families of matrices
and frames that satisfy the uncertainty principle.
These include: random orthogonal matrices, random partial Fourier matrices,
and a large family of matrices with independent entries (subgaussian
matrices), in particular random Gaussian and Bernoulli matrices.

\section{Kashin's representations}          \label{s: Kashin}

\subsection{Frame representations.} A sequence $(u_i)_{i=1}^N\subset \C^n$ is called
a tight frame if it satisfies Parseval's identity
\begin{equation}                                          \label{frame def}
\|x\|_2^2 =\sum_{i=1}^N |\< x,u_i \> | ^2
 \text{ for all $x\in \C^n$.}
\end{equation}
This definition differs by a constant normalization factor from  one which is
often used in the literature, but \eqref{frame def} will be more convenient
for us to work with.

A frame $(u_i)_{i=1}^N\subset \C^n$
can be identified with the $n \times N$ {\em frame matrix} $U$ whose columns are $u_i$.
The following properties are easily seen to be equivalent:
\begin{enumerate}
 \item $(u_i)_{i=1}^N$ is a tight frame for $\C^n$;
 \item every vector $x \in \C^n$ admits frame representation \eqref{frame};
 \item the rows of the frame matrix $U$ are orthonormal;
 \item \label{frame=projection}
   $u_i = P h_i$ for some orthonormal basis $(h_i)_{i=1}^N$ of $\C^N$,
   where $P$ is the orthogonal projection in $\C^N$ onto $\C^n$.
\end{enumerate}

When $N > n$, the tight frames are linearly dependent systems, so
various coefficients $a_i$ of the frame representation may carry common
 information about vector  $x\in \C^n$
This makes frames withstand noise
in coefficients better than orthonormal bases, see \cite{D, GVT, C}.
However, using  frame representation \eqref{frame}
may not always be the best way to use the frame redundancy.
Some coefficients $a_i$ may be much bigger than  others,
and thus carry more information about $x$.
In order to help information spread in the most uniform way,
one should try to make all coefficients of the same magnitude.
Such representations will be called Kashin's representations.

\subsection{Kashin's representations}
Consider a sequence $(u_i)_{i=1}^N\subset \C^n$.
We  say that the expansion
\begin{equation}              \label{Kashin's representation}
   x = \sum_{i=1}^N a_i u_i, \quad \max_i |a_i| \le \frac{K}{\sqrt{N}} \|x\|_2
\end{equation}
is a {\em Kashin's representation with level $K$} of  vector $x \in \C^n$.
\medskip

Kashin's representation produce the smallest possible dynamic range of the coefficients,
which is $\sqrt{n}$ smaller than the dynamic range of the frame representations.
This is the content of the following simple observation:

\begin{observation}[Optimality]                     \label{optimality}
 Let $(u_i)_{i=1}^N$ be a tight frame in $\C^n$. Then:

 (a) There exists a vector $x \in \C^n$ for which
 the coefficients $a_i = \< x, u_i\> $ of the frame representation \eqref{frame} satisfy
 $$
 \max_i |a_i| \ge  \sqrt{\frac{n}{N}} \|x\|_2.
 $$

 (b) For every vector $x \in \C^n$, every representation of the form
 $x = \sum_{i=1}^N a_i u_i$ satisfies
 $$
 \max_i |a_i| \ge \frac{1}{\sqrt{N}} \|x\|_2.
 $$
\end{observation}

\begin{proof}
(a) Since the tight frame satisfies $\sum_{i=1}^N \|u_i\|_2^2 = n$, one has
$\max_i \|u_i\|_2 \ge \sqrt{n/N}$. From this part (a) follows.

(b)  The correspondence  between tight frames and orthonormal bases 
(property 4) above) yields 
Bessel's inequality $\|x\|_2 \le ( \sum_{i=1}^N |a_i|^2 )^{1/2}$,
from which part (b) follows.
\end{proof}

Not every tight frame admits Kashin's representations with constant level $K$;
this is clear if one considers an orthonormal basis in $\C^n$ repeated
$\sim N/n$ times and properly normalized.
Nevertheless, some natural classes of frames do have this property.

We start with the following existence result.

\begin{theorem}[Existence]     \label{Kashin exists}
 There exist tight frames in $\C^n$ with arbitrarily small redundancy $\l = N/n > 1$,
 and such that every vector $x \in \C^n$ has  Kashin's representation
 with level $K$ that depends on $\lambda$ only (not on $n$ or $N$).
\end{theorem}

\begin{proof} This statement is in essence a reformulation of the classical result
from  geometric functional analysis due to Kashin \cite{K 77}
(with an optimal dependence $K = O\big(\sqrt{\frac{1}{\l-1} \log \frac{1}{\l-1}}\big)$ given later by Garnaev and Gluskin \cite{GG}).
To see this, we shall look at Kashin's representations from the geometric viewpoint.
Let $Q^N = \{x : \|x\|_\infty \le 1\}$ and $B^n = \{x : \|x\|_2 \le 1\}$
stand for the unit cube and unit Euclidean ball in $\C^N$ and $\C^n$ respectively.
The  observation  below  follows directly from the definition of Kashin's representations.

\begin{observation}                \label{Kashin=projections}
 {\em (Kashin's representations and projections of the cube):}
 Consider a tight frame  $(u_i)_{i=1}^N$ in $\C^n$ and a number $K > 0$.
 The following are equivalent:

 (i) Every vector $x \in \C^n$ has a Kashin's representation of level $K$
 with respect to the system $(u_i)_{i=1}^N$;

 (ii) The $n \times N$ matrix $U$ whose columns are $u_i$ satisfies
 \begin{equation}                               \label{ball in cube}
   B^n \subseteq \frac{K}{\sqrt{N}} U(Q^N).
 \end{equation}
\end{observation}

Inclusion \eqref{ball in cube} yields an equivalence
\begin{equation}                \label{Kashin geometric}
B^n \subseteq \frac{K}{\sqrt{N}} U(Q^N) \subseteq K B^n,
\end{equation}
the second inclusion holds trivially.
Since the rows of the frame matrix $U$ are orthonormal, the operator $U:\C^N \to \C^n$
is unitarily equivalent to an orthogonal projection. We thus may say that $U$
{\em realizes Euclidean projection of the cube.}
We refer the reader to  \cite{Pi} Section~6 for more thorough discussion of this topic.

Kashin's theorem \cite{K 77} states that there exists an orthogonal projection
of the unit cube in $\C^N$ onto a subspace of dimension $n$, which is equivalent
to a Euclidean ball and the coefficient $K$ depends on the redundancy $\l = N/n$ only.
In other words, there exists an $n \times N$ matrix $U$ whose rows are orthonormal and which satisfies
\eqref{Kashin geometric}.

The first inclusion in \eqref{Kashin geometric}
means that the columns $u_i$ of the matrix $U$
form a system for which every vector has a Kashin's representation.
Since the rows of $U$ are orthonormal, $(u_i)$ is a tight frame.
This proves Theorem~\ref{Kashin exists}.
\end{proof}

\medskip

In geometric functional analysis, many classes of matrices $U$
are known to realize Euclidean projections of the cube as
in \eqref{Kashin geometric}.
We discuss them in more details in Section~\ref{s: known examples}.
In fact we will see that   random matrix $U$ with orthonormal rows picked with respect to a rotationally
invariant distribution satisfies \eqref{Kashin geometric} with high probability.

\medskip

\noindent {\bf Remark} Since the level $K$ of Kashin's representation depends on redundancy only,
this representation become especially efficient in high dimensions when when the factor $\sqrt{n}$ in the expression 
for the dynamic range of the frame expansion overpowers the value of $K$ (which ideally is 
$O\left ( \sqrt {\frac 1 {\lambda-1}\log \frac 1 {\lambda-1}}\right )$). Therefore we are interested mainly in low redundant frames
just in order to avoid getting too large volumes of information to be transmitted.

\subsection{Stability, vector quantization}            \label{s: vector quantization}
Kashin's representations have   maximal power to reduce errors in the coefficients.
Indeed, consider a tight frame $(u_i)_{i=1}^N$ in $\C^n$, but instead of using
frame representations we shall use Kashin's representations with
some constant level $K = O(1)$.
So we represent a vector $x \in \C^n$, $\|x\|_2 \le 1$, with its Kashin's coefficients
$(a_1,\ldots,a_N) \in \C^N$, $|a_i| \le K/\sqrt{N}$. Assume these coefficients
are damaged (due to quantization, losses, flips of bits, etc.) and we only know
noisy coefficients $(\hat{a}_1,\ldots,\hat{a}_N) \in \C^N$.
When we try to reconstruct  $x$ from these  coefficients as
$\hat{x} = \sum_{i=1}^N \hat{a}_i u_i$, the accuracy of this
reconstruction is
\begin{IEEEeqnarray}{c}            \label{quantization error does not grow}
\| x - \hat{x}\|_2
=  \Big\| \sum_{i=1}^N (a_i - \hat{a}_i) u_i \Big\|_2    \notag \\
\le
\Big( \sum_{i=1}^N |a_i - \hat{a}_i|^2 \Big)^{1/2}.
\end{IEEEeqnarray}

Combined with the fact that the coefficients $a_i$ have the dynamic range $O(1/\sqrt{N})$,
this yields  greater robustness of Kashin's representations with respect to noise,
and in particular to {\em quantization errors}.
Suppose we need to quantize a vector $x \in \C^n$. We may do this by quantizing each coefficient
$a_i$ separately by performing a uniform scalar quantization
of the dynamic range $[-K/\sqrt{N}, K/\sqrt{N}]$ with, say, $L$ levels.
The quantization error for each coefficient is thus $|a_i - \hat{a_i}| \le K / L \sqrt{N}$.
By \eqref{quantization error does not grow}, this produces the overall quantization error
$$
\|x - \hat{x}\|_2 \le K/L = O(1/L).
$$
Similar quantization of frame representations \eqref{frame} would only give the bound
$$
\|x - \hat{x}\|_2 \le \sqrt{n}/L
$$
because its dynamic range is $\sqrt{n}$ larger than that of Kashin's representations
(by Observation~\ref{optimality}).

Kashin's decompositions also withstand {\em arbitrary errors} made to
a small fraction of the coefficients $a_i$. These may include losses of coefficients
and arbitrary flips of bits. Suppose at most $\d N$ coefficients
$(a_1, \ldots, a_N)$ are damaged in an arbitrary way, which results
in coefficients $(\hat{a}_1,\ldots,\hat{a}_N)$. Since all $|a_i| \le K/\sqrt{N}$,
we can assume (by truncation) that all $|\hat{a_i}| \le K/\sqrt{N}$. When we  
reconstruct  $x$ from these damaged coefficients (as before),
the accuracy of this reconstruction can be estimated using \eqref{quantization error does not grow}
as
$$
\| x - \hat{x}\|_2  \le  2 K \sqrt{\d} = O(\sqrt{\d}).
$$
Thus the reconstruction error is small whenever the (related) number of damaged coefficients $\d$ is small.

By Theorem~\ref{Kashin exists}, the maximal error reduction effect is achieved
using frames with only a {\em constant redundancy}, in fact  any redundancy
factor $\l = N/n > 1$ has the error reduction power
of maximal possible order. This is in contrast with traditional methods, in which increasing redundancy of
the frame gradually reduces the  representation error.

\section{Computing Kashin's representations}        \label{s: computing Kashin}

Computing the coefficients $a_i$ of Kashin's representation \eqref{Kashin's representation}
of a given vector $x$ can be described as a linear feasibility problem,
which can be solved in (weakly) polynomial time using linear
programming methods.

In this paper, we   take a different approach to computing Kashin's representations,
by establishing their connection with the uncertainty principle.
This will have several advantages over the linear programming approach:

\begin{enumerate}
 \item Whenever a frame $(u_i)$ satisfies the uncertainty principle,
one can effectively transform every frame representations into
Kashin's representation.
This will take $O(\log N)$ multiplications of the matrix $U$ by a vector.

 \item The uncertainty principle will thus be a guarantee that a given frame $(u_i)$
yields  Kashin's representation for every vector.
This can help to identify frames that yield Kashin's representations.

 \item The algorithm to transform frame representations into Kashin's representations
is simple, natural, and robust. It has a potential to be implemented on analog devices.
Followed by some robust scalar quantization of coefficients (such as
one-bit $\beta$-quantization \cite{DDGV1, DDGV2}), this algorithm may be used for
robust one-bit vector quantization schemes for analog-to-digital conversion.
\end{enumerate}

\subsection{The uncertainty principle}

The classical uncertainty principle says that a function and its Fourier transform
cannot be simultaneously well-localized. We refer the reader to fundamental monograph
\cite{HJ}  for history survey and also for numerous realization of this heuristic rule.
In particular a variant of the uncertainty principle due to Donoho and Stark
\cite{DS} states that if $f \in L_2(\R)$  is "almost concentrated" on a measurable set $T$
while its Fourier transform $\hat{f}$  is "almost concentrated" on a measurable set
$\Omega$, then then the product of measures $|T| |\Omega|$ admits a natural low
bound.  Donoho and Stark proposed applications of this
principle for signal recovery \cite{DS}.

For signals on discrete domains  no satisfactory version of the
uncertainty principle was known until recently.
For the discrete Fourier transform in $\C^N$
the uncertainty principle    states that $|\supp(x)| |\supp(\hat{x})| \ge N$
for all $x \in \C^N$ (see \cite{DS}).
This inequality is  sharp --
both terms in this product can be of order $\sqrt{N}$.

In papers by Candes, Romberg and Tao \cite{CRT, CT, CR}
and by Rudelson and Vershynin \cite{RV, RV full},
a much stronger discrete uncertainty principle was established
for {\em random} sets of size proportional to $N$.
Moreover, one of these sets  (say support of the signal in frequency domain)
can be arbitrary
(non-random), and the other (random  support in time domain) can be almost the whole domain.
The following result is a consequence of \cite{RV, RV full}:

\begin{theorem}[uncertainty principle]              \label{UP Fourier intro}
 Let $N = (1+\mu) n$ for some integer $n$ and $\mu \in (0,1)$.
 Consider a random subset $\Omega$ of $\{0,\ldots,N-1\}$ of average
 cardinality $n$, which is obtained from independent random $\{0,1\}$-valued
 variables $\d_0,\ldots \d_{N-1}$ with $\E \d_i = n/N$ as $\Omega := \{i: \; \d_i =1\}$.
 Then $\Omega$ satisfies the following with high probability.
 For every $z \in \C^N$,
 $$
 \supp(z) \subseteq \Omega  \ \ \ \text{implies} \ \ \
 |\supp(\hat{z})| > \d N,
 $$
 where $\d = c \mu^2 / \log^2 N$ and $c>0$ is an absolute constant.

 Moreover, for every $x \in \C^N$, $|\supp(x)| \le \d N$, one has
 \begin{equation}                      \label{x on Omega}
 \|\hat{x} \cdot \one_\Omega\|_2 \le (1-c\mu) \|x\|_2,
 \end{equation}
 where $\one_\Omega$ denotes the indicator function of $\Omega$.
\end{theorem}

The first, qualitative, part of the theorem easily follows
from the second, quantitative part with $z = \hat{x}$.
If $\supp(\hat{x}) \subseteq \Omega$
and $|\supp(x)| \le \d N$ then, by the second part,
$\|\hat{x}\|_2 = \|\hat{x} \cdot \one_\Omega\|_2 < \|x\|_2$,
which would contradict Parseval equality.

We can regard inequality \eqref{x on Omega} as a property
of the partial Fourier matrix $U$,
which consists of the rows of the DFT  (discrete Fourier transform)
matrix $\Phi$ indexed by
the random set $\Omega$. Then \eqref{x on Omega} says that
$\|U x\|_2 \le (1-c\mu) \|x\|_2$ for all vectors
$x \in \C^N$ such that $|\supp(x)| \le \d N$.

Now we can abstract from the harmonic analysis in question and introduce a
general uncertainty principle (UP) as a property of matrices.

\begin{definition}[UP for matrices]      \label{UP for matrices}
 An $n \times N$ matrix $U$ satisfies the uncertainty principle with
 parameters $\eta, \d \in (0,1)$ if, for $x \in \C^N$,
 \begin{equation}       \label{Eq: UP frames}
 \ |\supp(x)| \le \d N \text{ implies } \|Ux\|_2 \le \eta \|x\|_2.
 \end{equation}
\end{definition}

We will only use the uncertainty principle for matrices $U$ with orthonormal
(or almost orthonormal)  rows, in which case it is always a nontrivial property.

\medskip

A related uniform uncertainty principle (UUP) was introduced by Candes and Tao
in the context of the sparse recovery problems \cite{CT2}.
The UUP with parameters $\e, \d \in (0,1)$ states that there exists $\l > 0$ such that,
for $x \in \C^N$, the condition $|\supp(x)| \le \d N$ implies
$$
\l(1-\e) \|x\|_2 \le \|Ux\|_2 \le \l(1+\e) \|x\|_2.
$$
See also \cite{CRT2, CRTV} for more refined versions.
Known also as the Restricted Isometry Condition, UUP was shown in \cite{CT2}
to be a guarantee that one can efficiently solve underdetermined systems
of linear equations $Ux = b$
under the assumption that the solution is sparse, $|\supp(x)| \le \d N$.
This is a part of the fast developing area of Compressed Sensing \cite{CS webpage}.

The uncertainty principle is a weaker assumption (thus easier to verify) than the UUP:

\begin{observation}
 For matrices with orthonormal rows,
 the UUP with parameters $\e, \d$ implies the uncertainty principle
 with parameters $\eta = \frac{1+\e}{1-\e} \sqrt{\frac{n}{N}}$, $\d$.
\end{observation}

\begin{proof}
Since the columns $u_i$ of the matrix $U$ satisfy $\sum_{i=1}^N \|u_i\|_2^2 = n$,
there exists a column with norm $\|u_i\|_2 \le \sqrt{n/N}$. This column is  a
preimage of some $1$-sparse unit vector $e_i$, i.e. $u_i = U e_i$ where
$e_i = (0,\ldots,0,1,\ldots,0)$ with $1$ on the $i$-th place.
Using the UUP for $x=e_i$ we obtain
$$
\l (1-\e) \le \|u_i\|_2 \le \sqrt{\frac{n}{N}}.
$$
Hence $\l \le \frac{1}{1-\e} \sqrt{\frac{n}{N}}$.
In view of this estimate, the upper bound in the UUP reads as follows:
$|\supp(x)| \le \d N$ implies
$$
\|Ux\|_2 \le \frac{1+\e}{1-\e} \sqrt{\frac{n}{N}} \|x\|_2.
$$
This is what we wanted to prove.
\end{proof}

The uncertainty principle can be reformulated as a property of
systems of vectors $(u_i)_{i=1}^N$, which form the columns of the matrix $U$.
We will use it for tight (or almost tight) frames, in which case it is
a nontrivial property:

\begin{definition}[UP for frames]          \label{UP frames}
 A system of vectors $(u_i)_{i=1}^N$ in $\C^n$ satisfies the uncertainty
 principle with parameters $\eta,\d \in (0,1)$ if
 \begin{equation}                                          \label{Hilbertian}
   \Big\| \sum_{i\in \Omega} a_i u_i \Big\|_2
   \le \eta \Big( \sum_{i\in \Omega} |a_i|^2 \Big)^{1/2}
 \end{equation}
 for every subset $\Omega \subset \{1,2,\ldots \, ,N\}$, $|\Omega|\leq \d N$.
\end{definition}

\subsection{Converting frame representations into Kashin's representations}
\label{s: frame to Kashin}

For every tight frame that satisfies the uncertainty principle,
one can convert frame representations into Kashin's representations.

The conversion  procedure is natural and fast.
We truncate the coefficients of the frame representation
\eqref{frame} of $x$ at level $M = \|x\|_2/\sqrt{\d N}$
in hope to achieve a Kashin's representation
with level $K = 1/\sqrt{\d}$. However, the truncated representation may sum up
to a vector $x^{(1)}$ different from $x$.
So we consider the residual $x - x^{(1)}$,
compute its frame representation and again truncate its coefficients,
now at a lower level $\eta M$. We continue this process of
expansion, truncation and reconstruction, each time reducing the truncation
level by the factor of $\eta$.

Using the uncertainty principle, we will be able to show that
the norm of the residual reduces by the factor of $\eta$ at each iteration.
So we can compute Kashin's representations of level $K = K(\eta,\d)$
with accuracy $\e$ in $O(\log(1/\e))$ iterations.
Our analysis of this algorithm will yield:

\begin{theorem}[frame to Kashin conversion]      \label{algorithm valid}
 Let $(u_i)_{i=1}^N$ be a tight frame in $\C^n$ which satisfies the Uncertainty
 Principle with parameters $\eta,\d$. Then each vector
 $x\in \C^n $ admits a Kashin representation of level $K=(1-\eta)^{-1}\d^{-1/2}$.
\end{theorem}

In order to prove this result, we introduce and study the truncation operator
for frame representations.
Given a number $M>0$, the one-dimensional truncation at level $M$ is defined
for $z \in \C \setminus \{0\}$ as
\begin{equation}                                    \label{one d truncation}
t_M(z) = \frac{z}{|z|}\min\{|z|, M\},
\end{equation}
and $t_M(0) = 0$.

Consider a frame $(u_i)_{i=1}^N$ satisfying the assumptions of the theorem.
For every $x \in \C^n$, we consider the frame representation
$$
x = \sum_{i=1}^N b_i u_i \text{ where } b_i = \< x, u_i \>
$$
and define the truncation operator on $\C^n$ as
\begin{multline}
  Tx=\sum_{i=1}^N \hat{b}_i u_i \text{ where } \hat{b}_i=t_M(b_i) \\
  \text{ and } M = \|x\|_2/\sqrt{\d N}.
\end{multline}

The uncertainty principle helps us to bound the residual of the truncation:

\begin{lemma}[Truncation]              \label{truncation error lemma}
 In the above notations, for every vector $x \in \C^n$ we have
 \begin{equation}
   \|x-Tx\|_2 \leq \eta \|x\|_2.
 \end{equation}
\end{lemma}

\begin{proof}
Let $x\in \C^n$. Consider the subset $\Omega \subseteq \{1,\ldots,N\}$
defined as
\[
\Omega = \{i : \; b_i\neq \hat{b}_i\} = \{i : \; |b_i| > M\}.
\]
By the definition of tight frame, we have
$$
\|x\|_2^2 = \sum_{i=1}^N |b_i|^2 > |\Omega| M^2,
$$
thus
$$
|\Omega| \le \|x\|_2^2 / M^2 = \delta N.
$$
Using the uncertainty principle, we can estimate the norm of the residual
\[
x - T_M x=\sum_{i\in \Omega}(b_i-\hat{b}_i)u_i
\]
as
\begin{multline*}
\|x - T_M x\|_2 \leq \eta \Big (\sum_{i\in \Omega}|b_i-\hat{b}_i|^2\Big )^{1/2}
 \le
 \\
 \eta \Big (\sum_{i\in \Omega}|b_i|^2\Big )^{1/2}
 \le \eta \Big (\sum_{i=1}^N|b_i|^2\Big )^{1/2}
 = \eta \|x\|_2.
\end{multline*}
This completes the proof.
\end{proof}
\medskip

\noindent{\bf Proof of Theorem 3.5}
Given $x\in \C^n$, for $k=1,2,\ldots$ we define the vectors
$$
x^{(0)}:=x, \quad x^{(k)}:=x^{(k-1)} - Tx^{(k-1)}.
$$
Then, for each $r=0,1,2,\ldots$ we have
\[
x = \sum_{k=0}^r Tx^{(k)} + x^{(r+1)}.
\]
It follows from Lemma \ref{truncation error lemma} by induction that
$\|x^{(k)}\|_2\leq \eta^k \|x\|_2$, thus
\[
x = \sum_{k=0}^\infty Tx^{(k)}.
\]
Furthermore, by the definition of the truncation operator $T$,
each vector $Tx^{(k)}$ has an expansion in the system $(u_i)_{i=1}^N$
with coefficients bounded by
$\|x^{(k)}\|_2/\sqrt{\d N} \le \eta^k \|x\|_2/\sqrt{\d N}$.
Summing up these expansions for $k=0,1,2,\ldots$,
we obtain an expansion of $x$
with coefficients bounded by
$(1-\eta)^{-1} \|x\|_2/\sqrt{\d N}$.
In other words, $x$ admits Kashin's representation with level
$K = (1-\eta)^{-1} \d^{-1/2}$.
This completes the proof.
\endproof

\medskip

The proof yields an algorithm to compute Kashin's representations:

\medskip

\textsc{Algorithm to compute Kashin's representations}
\nopagebreak \\
  {\bf Input:}
     \begin{itemize}
       \item A tight frame $(u_i)_{i=1}^N$ in $\C^n$
         which satisfies the uncertainty principle with parameters $\eta,\d \in (0,1)$.
       \item A vector $x \in \C^n$ and a number of iterations $r$.
     \end{itemize}
  {\bf Output:}
    Kashin's decomposition of $x$ with level $K= (1-\eta)^{-1}\d^{-1/2}$
    and with accuracy $\eta^r \|x\|_2$.
     Namely, the algorithm finds coefficients $a_1, \ldots, a_N$
     such that
     \begin{multline}
     \Big\|x - \sum_{i=1}^N a_i u_i \Big\|_2  \le \eta^r \|x\|_2, \\
     \max_i |a_i| \le \frac{K}{\sqrt{N}} \|x\|_2.
     \end{multline}
 \smallskip
 Initialize the coefficients and the truncation level:
      $$
      a_i \leftarrow 0, \  i=1,\ldots,N;
      \quad M \leftarrow \frac {\|x\|_2}{\sqrt{\d N}}.
      $$
 Repeat the following $r$ times:
     \begin{itemize}
       \item Compute the frame representation of $x$ and truncate at level $M$:
         $$
         b_i \leftarrow \< x,u_i \> , \;
           \hat{b}_i \leftarrow t_M(b_i), \;  i=1,\ldots,N.
         $$
       \item Reconstruct and compute the error:
         $$
         Tx \leftarrow \sum_{i=1}^N \hat{b}_i u_i; \quad
         x \leftarrow x-Tx.
         $$
       \item Update Kashin's coefficients and the truncation level:
         \begin{align*}
         &a_i \leftarrow a_i + \sqrt{N} \, \hat{b}_i, \  i=1,\ldots,N; \\
         &M \leftarrow \eta M. \qquad \qquad
         \end{align*}
     \end{itemize}

\bigskip

\begin{remark} {\em (Redistributing information). }
One can view this algorithm as a method of redistributing information
among the coefficients. At each iteration, it ``shaves off''
excessive information from the few biggest coefficients (using truncation)
and redistributes this excess more evenly.
This process is continued until all coefficients have a fair share
of the information.
\end{remark}

\medskip

\begin{remark}{\em (Computing exact Kashin's representations).}
With a minor modification, this algorithm
can compute an {\em exact} Kashin's representation after $r = O(\log N)$ iterations.
We just do not need to truncate the coefficients $b_i$ during the last iteration.

Indeed, for such $r$, the error factor satisfies $\eta^r \le \frac{K}{\sqrt{N}}$.
Thus, during $r$-th iteration the frame coefficients $b_i$
are all bounded by $\frac{K}{\sqrt{N}} \|x\|_2$, where $x$ is the initial input vector.
So $b_i$ are already sufficiently small,
and we will not apply the truncation at the last iteration.
This yields an exact Kashin's representation of $x$ with $K' = 2K$.
\end{remark}

\medskip

\begin{remark}{\em (Robustness).}
 The algorithm above is robust in the sense of \cite{DD}.
\end{remark}

Specifically, the truncation operation \eqref{one d truncation}
may be impossible to realize  on a physical signal exactly,
because it is expensive to build an analog scheme
that has an exact phase transition at the truncation level $|z| = M$.
A robust algorithm should not rely on any assumptions on exact phase
transitions of the operations it uses.
Scalar quantizers that are robust in this sense were
first constructed by Daubechies and DeVore in \cite{DD} and
further developed in \cite{G, DDGV1, DDGV2}.

Our algorithm is also robust in the following sense: the exact truncation
$t_M$ can be replaced by any approximate truncation.
Such an approximate truncation at level $1$ can be any function
$t(z): \C \to \C$ which satisfies for some $\nu, \tau \in (0,1)$:
\begin{eqnarray}                                   \label{truncation}
 |z - t(z)| \le
   \begin{cases}
     \nu |z| & \text{if $|z| \le \tau$}, \\
     |z|     & \text{for all $z$,}
   \end{cases} \\
  \text{and $|t(z)| \le 1$ for all $z$}.    \nonumber
  \end{eqnarray}
The approximate truncation at level $M$ is defined as
$t_M (u) := M \, t(\frac{u}{M})$.
An analysis similar to that above yields:

\begin{theorem}  {\em (Approximate truncation)}
 The above algorithm  remains valid if one replaces exact truncation
 by any approximate one and also adjust the parameters:
level $M$ should be replaced with $M' = \t^{-1} M$,
 parameter $\eta$ should be replaced with $\eta' = \sqrt{\eta^2 + \nu^2}$,
 finally,  level $K$  is replaced with $K' = \t^{-1} (1-\eta')^{-1} \d^{-1/2}$,
 provided that $\eta' < 1$.
\end{theorem}

Moreover, the approximate truncation can be different each  time it is
called by the algorithm, provided that it satisfies  \eqref{truncation}.
This facilitates   the algorithm implementation on analog devices.
In particular, one can use this algorithm
to build robust vector quantizers for analog-to-digital conversion.

\bigskip

\begin{remark}{\em (Almost tight frames).}
 Similar results also hold for frames that are almost, but not exactly, tight.
 This is important for natural classes of frames, such as random gaussian
 and subgaussian frames (see Theorem~\ref{subgaussian}).
\end{remark}

\begin{definition}
 For $\e \in (0,1)$, a sequence $(u_i)_{i=1}^N\subset \C^n$ is
 called an $\e$-tight frame if
 \begin{multline}   \label{e tight frame}
 (1-\e)\|x\|_2\leq \Big( \sum_{i=1}^N |\< x,u_i \> | ^2 \Big)^{1/2}  \\
           \leq (1+\e)\|x\|_2 \quad \text{for all $x\in \C^n$.}
 \end{multline}
\end{definition}

\smallskip

An analysis similar to that above yields:

\begin{theorem}                     \label{almost tight}
 Let  $(u_i)_{i=1}^N\subset \C^n$ be an $\e$-tight frame,
 which satisfies the uncertainty principle with parameters $\eta$ and
 $\delta$.
 Then Theorem~\ref{algorithm valid} and the algorithm above
 are valid for $M$ replaced with $M' = \sqrt{1+\e}\,M$
 and $\eta$ replaced with $\eta' = \sqrt{1+\e}\;\eta + \e$,
 provided that $\eta' < 1$.
\end{theorem}

\begin{remark}{\em (History).}
The idea behind Theorem~\ref{algorithm valid}
is certainly not new. Gluskin \cite{G unpublished} suggested to use properties
that involved only $\|\cdot\|_2$ norms (like our uncertainty principle)
to deduce results on Euclidean sections of $\ell_1^n$ (which by duality
is equivalent to Euclidean projections \eqref{ball in cube} of a cube).
A similar idea was essentially used by Talagrand in his work on the
$\Lambda_1$ problem \cite{T}.

The algorithm to compute Kashin's representations resembles the
Chaining Algorithm of \cite{GSTV}, which also detects a few
biggest coefficients and iterates on the residual, but it serves
to {\em find} all big coefficients rather than to spread them out.
\end{remark}


\section{Matrices and frames that satisfy the uncertainty principle}
\label{s: matrices UP}

In this section, we give examples of matrices (equivalently, frames)
that satisfy the uncertainty principle.
By Observation~\ref{Kashin=projections},
such $n \times N$ matrices $U$ realize Euclidean projection of the
cube \eqref{ball in cube}. Equivalently, these frames $(u_i)_{i=1}^N$
(the columns of $U$) yield quickly computable Kashin's representations
for every vector $x \in \C^n$.


\subsection{Matrices known to realize Euclidean projections of the cube}
\label{s: known examples}
Much attention has been paid to   Euclidean projections of the cube
\eqref{ball in cube} in geometric functional analysis.
Results in the literature are usually stated in the dual form,
about $n$-dimensional Euclidean subspaces of $\ell_1^N$.

Kashin proved \eqref{ball in cube} for  random orthogonal $n \times N$
matrix $U$ (formed by the first $n$ rows of a random matrix in $\OO(N)$),
with $N = \l n$ for arbitrary $\l > 1$, and
with exponential probability (\cite{K 77}, see also \cite{Pi} Section 6.)
The level $K$ \eqref{ball in cube} depends only on $\l$;
an optimal dependence was given later by Garnaev and Gluskin \cite{GG}).

A similar result holds for $U = \frac{1}{\sqrt{N}} \Phi$, where
$\Phi$ is a random Bernoulli matrix, which means that the entries
of $\Phi$ are $\pm 1$ symmetric independent random variables.
Schechtman \cite{S} first proved this   with $N = O(n)$, and in
\cite{LPRTV} this result is generalized for $N = \l n$ with arbitrary $\l > 1$.
The dependence  $K$ on $\l$ was improved recently in \cite{AFMS}.
In fact, these results hold for a quite general class of subgaussian matrices
(which includes Bernoulli and Gaussian random matrices).

In is unknown whether Kashin's theorem holds for  partial Fourier
matrix; this conjecture is known as the $\Lambda_1$ problem.
Consider the Discrete Fourier Transform in $\C^{N}$, where $N =
O(n)$, given by the orthogonal $N \times N$ matrix $\Phi$. It is
unknown whether there exists a submatrix $U$ which consists of
some $n$ rows of $\Phi$ and such that it realizes an Euclidean
projection of the cube in the sense of \eqref{ball in cube}.

In the positive direction, a partial result due to Bourgain,
later reproved by Talagrand with a general method \cite{T},
states that a random partial Fourier matrix $U$ satisfies \eqref{Kashin geometric}
with high probability for $N = O(n)$ and $K = O(\sqrt{\log(N) \log \log(N)})$.
It was recently proved in \cite{GMPT}
that Bourgain's result holds for arbitrarily small redundancy, that is for
$N = \l n$ with arbitrary $\l > 1$,
however at the cost of a slightly worse logarithmic factor in $K$.
A similar result can also be deduced from Theorem~\ref{Fourier} below
(along with Theorem~\ref{algorithm valid} and ~\ref{Kashin=projections}),
which is a consequence of the uncertainty principle in \cite{RV, RV full}.

No explicit constructions of matrices $U$ are known.
However, there exists small space constructions that
use a small number of random bits \cite{AM, I, I2}.

\subsection{Random orthogonal matrices}         \label{s: orthogonal}

We consider random $n\times N$ matrices whose rows are orthonormal.
Such matrices can be obtained by selecting the first $n$ rows of orthogonal
$N\times N $ matrices. Indeed, denote by $\OO(N)$ the space of all
orthogonal $N\times N$ matrices with the normalized Haar measure. Then
\begin{equation}                                           \label{matrix projection}
\OO(n\times N)=\{ P_nV; V\in \OO(N)\},
\end{equation}
where $P_n:\C^N\to \C^n$ is the orthogonal projection on the first $n$ coordinates.
The probability measure on $\OO(n\times N)$ is induced by the Haar
measure on $\OO(N)$.

\begin{theorem} {\em (UP for random orthogonal matrices)}  \label{random orthogonal}

 Let $\mu>0$ and $N = (1+\mu) n$.
 Then, with probability at least $1-2\exp(-c \mu^2 n)$,
 a random orthogonal $n \times N$ matrix $U$
 satisfies the uncertainty principle with the  parameters
 \begin{equation}               \label{random frame parameters}
 \eta = 1 - \frac{\mu}{4}, \ \ \
 \d = \frac{c\mu^2}{\log(1/\mu)},
 \end{equation}
 where $c > 0$ is an absolute constant.
\end{theorem}

\medskip

\noindent {\bf Remark.} Assumption $\mu>0$ is not essential; just expressions for
$\eta$ and $\d$ will look differently. We are most interested in small values of  $\mu$
when redundancy is small.

\medskip

The proof of Theorem~\ref{random orthogonal} uses a standard scheme
in geometric functional analysis -- the concentration inequality on
the sphere followed by an $\e$-net argument.
Denote by $S^{N-1}$ and $\s_{N-1}$ the unit Euclidean sphere in $\C^N$ and
the normalized Lebesgue measure on $S^{N-1}$.

\begin{lemma}                            \label{U tail}
 For arbitrary $t > 0$, $x\in S^{N-1}$, we have
 $$
 \P \left \{ \|Ux\|_2 > (1+t) \sqrt{\frac{n}{N}}\right \} \le 2 \exp(-c_1 t^2 n),
 $$
 where $c_1 > 0$ is an absolute constant.
\end{lemma}

\begin{proof}
We use  representation \eqref{matrix projection} and also the fact that
$z = Vx$ is a random vector uniformly distributed on  $S^{N-1}$.
Thus $Ux$ is distributed identically with $P_n z$.
We also have
\begin{multline*}
E:=\int_{S^{N-1}}\| P_n z\|_2 \, d \s_{N-1}(z)
\le
\\
\Big(\int_{S^{N-1}} \| P_n z\|_2^2 \; d \s_{N-1}(z) \Big)^{1/2}
= \sqrt{\frac{n}{N}}.
\end{multline*}
The map $z \mapsto \|P_nz\|$ is a $1$-Lipschitz function on $S^{N-1}$.
The  concentration inequality (see e.g. \cite{L} Section~1.3) then implies
that this function is well concentrated about its average value $E$:
\begin{multline*}
\P \{ \|Ux\|_2 > E+ u \}
\leq \\
\s_{N-1} (z \in S^{N-1} :\; |\|P_nz \|_2- E| > u \})
\le
\\
2 \exp(-c u^2 N).
\end{multline*}
Choosing $u = t \sqrt{n/N}$ completes the proof.
\end{proof}

\noindent {\bf Proof of Theorem~\ref{random orthogonal}. }
Assume that  $\eta$ and $\d$ satisfy the assumptions \eqref{random frame parameters}.
We have to prove that  \eqref{Eq: UP frames} holds with probability
at least $1-\exp(-c \mu^2 n)$ .

Consider the set
$$
S := \{ x \in S^{N-1}, \; |\supp(x)| \le \d N \} .
$$
We have
$$
S=\bigcup_{|I| \le \d N} S_I,
$$
here the union is  taken over all subsets $I$ of $\{1,\ldots,N\}$ of cardinality
at most $\d N$,
and $S_I = S^{N-1} \cap \C^I$ is the set of all unit vectors whose supports
lie in $I$.
Let $\e > 0$. For each $I$, we can find an $\e$-net of $S_I$ in the Euclidean norm,
and of cardinality at most $(3/\e)^{\d N}$ (see e.g. \cite{Pi} Lemma~4.16).
Taking the union over all sets $I$ with $|I| = \lceil \d N \rceil$,
we conclude by the Stirling's bound on the binomial coefficients
that there exists an $\e$-net $\NN$ of $S$ of cardinality
$$
|\NN| \le \binom{N}{\lceil \d N \rceil} \Big(\frac{3}{\e}\Big)^{\d N}
\le \Big( \frac{3e}{\e\d} \Big)^{\d N}.
$$
Then using Lemma~\ref{U tail}, we obtain
\begin{multline*}
\P \{ \exists y \in \NN :\; \|Uy\|_2 > (1+t) \sqrt{\frac{n}{N}} \}
\le
\\
|\NN| \cdot 2 \exp(-c_1 t^2 n).
\end{multline*}
Every $x \in S$ can be approximated by some $y \in \NN$ within $\e$
in the Euclidean norm, and since $U$ has norm one, we have
$$
\|Ux\|_2 \le \|Uy\|_2 + \|U(x-y)\|_2 \le \|Uy\|_2 + \e.
$$
Therefore
\begin{multline}            \label{prob}
\P \{ \exists x \in S :\; \|Ux\|_2 > (1+t) \sqrt{\frac{n}{N}} + \e \}
\le
\\
|\NN| \cdot 2 \exp(-c_1 t^2 n).
\end{multline}
It now remains to choose parameters appropriately.
Let $t = \mu/5$ and $\e = \mu/8$.
Then since $N/n = 1+\mu$ and by the assumption on $\eta$ in
\eqref{random frame parameters}, we have
$$
(1+t) \sqrt{\frac{n}{N}} + \e \le \eta.
$$
Also, we can estimate the probability in \eqref{prob} as
\begin{multline}                    \label{probab}
|\NN| \cdot 2 \exp(-c t^2 n)
\le
\\
\Big(\frac{24e}{\d \mu} \Big)^{\d N}
 \cdot 2 \exp(-c_2 t^2 n) \le
 \\
2 \exp \big[ (2\d \log(24e/\d \mu) - c_2 \mu^2) n \big],
\end{multline}
where $c_2 = c_1/25$. By our choice of $\d$, the right hand side of
\eqref{probab} is bounded by $2\exp(-c \mu^2 n)$, where $c > 0$
is an absolute constant.
We conclude that
$$
\P \{ \exists x \in S :\; \|Ux\|_2 > \eta\} \le 2\exp(-c \mu^2 n).
$$
This completes the proof.
\endproof

\subsection{Random partial Fourier matrices}            \label{s: Fourier}

An important class of matrices that satisfy the uncertainty principle
can be obtained by selecting $n$ random rows of an arbitrary
orthogonal $N\times N$ matrix $\Phi$ whose entries are $O(N^{-1/2})$.
Here $n$ can be an arbitrarily big fraction of $N$, so the Uncertainty
Principle will hold for almost square random submatrices.
This class includes  random partial Fourier matricies,  multiplication
by such matrix corresponds to sampling $n$ random frequencies of a signal.

More precisely, we select rows of $\Phi$ using random selectors
$\d_1,\ldots,\d_N$ -- independent Bernoulli random variables,
which take value $1$ each with probability $n/N$.
The selected rows will be indexed by a random subset
$\Omega = \{i : \d_i = 1\}$ of $\{1,\ldots,N\}$,
whose average cardinality is $n$.

\begin{theorem}     \label{Fourier} {\em (UP for random partial Fourier matrices)}

 Let $\Phi$ be an orthogonal $N \times N$ matrix with uniformly
 bounded entries: $|\Phi_{ij}| \le \a N^{-1/2}$ for some constant $\a$
 and all $i,j$.
 Let $n$ be an integer such that $N = (1+\mu) n$ for some $\mu \in (0,1]$.
 Then for each $p \in (0,1)$ there exists a constant $c=c(p,\alpha) > 0$
 such that the following holds.

 Let $U$ be a submatrix of $\Phi$ formed by selecting a subset of
 the rows of average cardinality $n$.
 Then, with probability at least $1-p$, the matrix $U$
 satisfies the uncertainty principle with parameters
 $$
 \eta = 1 - \frac{\mu}{4}, \ \ \
 \d = \frac{c \mu^2}{\log^4 N}.
 $$
\end{theorem}

\medskip

Theorem~\ref{Fourier} is a direct consequence of a slightly
stronger result established in \cite{CT} and improved in \cite{RV, RV full}.
For an operator $U$ on a Euclidean space, $\|\cdot\|$ will denote its
operator norm.

\begin{theorem}  \label{UUP Fourier}
{\em (UUP for partial Fourier matrices \cite{RV, RV full})}

Assume the hypothesis of Theorem \ref{Fourier} is satisfied.
 Then there exists a constant $C = C(\a) > 0$ such that the following holds.
 Let $r > 0$ and $\e \in (0,1)$ be such that
 $$
 n \ge C \Big(\frac{r \log N}{\e^2}\Big)
       \log \Big(\frac{r \log N}{\e^2}\Big)
       \log^2 r.
 $$
 Then the random submatrix $U$ satisfies:
 \begin{equation}                  \label{E norm}
   \E \sup_{|T| \le r} \| \id_T - \frac{N}{n} U_T^* U_T \| \le \e.
 \end{equation}
 Here the supremum is taken over all subsets $T$ of $\{1,\ldots,N\}$
 with at most $r$ elements, $U_T$ denotes the submatrix of $U$
 that consists of the columns of $U$ indexed in $T$,
 and $\id_T$ denotes the identity on $\C^T$.
\end{theorem}
\medskip

\noindent {\bf Proof of Theorem~\ref{Fourier}.}
Observe that for a linear operator $A$ on $\C^N$ one has
\begin{multline}                    \label{duality}
\|\id - A^*A\|
= \sup_{x \in \C^N, \|x\|_2 = 1} \big| \< (\id-A^*A)x,x \> \big|
\\
= \sup_{x \in \C^N, \|x\|_2 = 1} \big| \|Ax\|_2^2 - \|x\|_2^2 \big|.
\end{multline}
We use this observation for $A = \sqrt{\frac{N}{n}} U_T$.
Since $U_T x = Ux$ whenever $\supp(x) \subseteq T$, we obtain
$$
\E \sup_{\substack{x \in \C^N, \,  \|x\|_2 = 1, \\ |\supp(x)| \le r }}
\Big| \frac{N}{n} \|Ux\|_2^2 - 1 \Big| \le \e.
$$
By Markov's inequality, with probability at least $1-p$
the random matrix $U$ satisfies:
\begin{multline*}
\ \ \ \ \ \ \Big| \frac{N}{n} \|Ux\|_2^2 - 1 \Big| < \e/p \\
\text{for all $x \in \C^N$, $|\supp(x)| \le r$, $\|x\|_2 = 1$.}
\end{multline*}
In particular, for such $U$, one has:
$
\|Ux\|_2 \le $
\begin{multline*}
\sqrt{1+\e/p} \sqrt{\frac{n}{N}} \, \|x\|_2
= \sqrt{\frac{1+\e/p}{1+\mu}} \, \|x\|_2
\\
\text{for all $x \in \C^N$, $|\supp(x)| \le r$.} \qquad
\end{multline*}
Then, if we set $\e = c p \mu$ for an appropriate absolute constant $c>0$,
we can bound the factor
$$
\sqrt{\frac{1+\e/p}{1+\mu}} \le 1 - \frac{\mu}{4} = \eta.
$$
This proves the uncertainty principle \eqref{Eq: UP frames} with $\d = r/N$.
To estimate $\d$ we note that the condition on $n$ in Theorem~\ref{UUP Fourier}
is satisfied if
\begin{equation}                    \label{r}
r \le \frac{c_1 \e^2 N}{\log^4 N}
\end{equation}
where $c_1 = c_1(\a) > 0$.
Since we have set $\e = c p \mu$, condition \eqref{r} is equivalent to
$$
\d \ge \frac{c \mu^2}{\log^4 N}
$$
where $c = c(\a) >0$.
This completes the proof of Theorem~\ref{Fourier}.
\endproof

\medskip

\begin{remarks}

1. {\em (Computing in almost linear time). }
The Fourier matrices can be used to compute Kashin's representations
in $\C^n$  in time almost linear in $n$. Indeed, let for example $N = 2n$.
The columns of the $n \times N$ partial Fourier matrix form a tight frame
in $\C^n$. By Theorem~\ref{Fourier} and Section~\ref{s: frame to Kashin},
we can convert a frame representation of every vector $x \in \C^n$
into a Kashin's representation with level $K = O(\log^2 n)$
 in time $O(n \log^2 n)$.
(Recall that the algorithm makes $O(\log n)$ multiplications by a partial
Fourier matrix, and each multiplication can be done using the fast
Fourier transform in time $O(n \log n)$).

\smallskip

2. The constant $c = c(\a,p)$ depends polynomially on $\a$
and polylogarithmically on $p$. The polynomial dependence on $\a$
is straightforward form the proof of Theorem~\ref{UUP Fourier}
in \cite{RV, RV full}. The proof above gives a polynomial dependence
on the probability $p$. To improve it to a polylogarithmic dependence,
one can use an exponential tail estimate, proved in \cite{RV full} Theorem~3.9,
instead of the expectation estimate \eqref{E norm}.

\smallskip

3. We stated Theorem~\ref{Fourier} in the range $\mu\in (0,1]$ which is
most interesting for us (where the redundancy factor is small).
A similar result holds for arbitrary $\mu > 0$.
\end{remarks}

\subsection{Subgaussian random matrices.}

A large family of matrices with independent random entries satisfies the uncertainty principle.

\begin{definition}
A random variable $\phi$ is called {\em subgaussian with parameter $\b$}
if
$$
\P \{ |\phi| > u \} \le \exp(1-u^2/\b^2) \ \ \
\text{for all $u > 0$}.
$$
\end{definition}
Examples of subgaussian random variables include Gaussian $N(0,1)$
random variables and bounded random variables.

\begin{theorem}    \label{subgaussian}
{\em (UP for random subgaussian matrices)}

 Let $\Phi$ be a $n \times N$ matrix whose entries are
 independent mean zero subgaussian random variables with parameter $\b$.
 Assume that $N = \l n$ for some $\l \ge 2$.
 Then, with probability at least $1 - \l^{-n}$,
 the random matrix $U = \frac{1}{\sqrt{N}} \Phi$
 satisfies the uncertainty principle with parameters
 \begin{equation}                  \label{subgaussian parameters}
 \eta = C\b \sqrt{\frac{\log \l}{\l}}, \ \ \   \d = \frac{c}{\l},
 \end{equation}
 where $C, c > 0$ are absolute constants.
\end{theorem}

\begin{remark}
 Theorem~\ref{subgaussian} and Lemma~\ref{almost orthogonal} below
 can be deduced from the recent works \cite{MPT, MPT short}. However, we feel
 that it would be helpful to include short and rather standard proofs
 of these results here.
\end{remark}

\medskip

Theorem~\ref{subgaussian} follows easily from an estimate on the operator
norm of subgaussian matrix.

\begin{lemma} (\cite{LPRT} Fact 2.4)         \label{Phi norm}
 Let $n \ge k$ and $\Phi$ be a $n \times k$ matrix whose entries are
 independent mean zero subgaussian random variables with parameter $\b$.
 Then
 \begin{equation}        \label{Eq:Phi norm}
 \P \{ \|\Phi\| > t \sqrt{n} \} \le \exp(-c_1 n t^2/\b^2)
   \end{equation}
 $ \text{for all $t \ge C_1 \b$}$,
here $C_1, c_1 > 0$ are absolute constants.
\end{lemma}

\noindent {\bf Proof of Theorem~\ref{subgaussian}. }
The uncertainty principle for the matrix $U$ with parameters $\eta, \d$
is equivalent to the following norm estimate:
$$
\sup_{|I| = \lceil \d N \rceil} \|\Phi_I\| \le \eta \sqrt{N},
$$
where the supremum is over all subsets $I \subset \{1,\ldots,N\}$
of cardinality $\lceil \d N \rceil$, and where $\Phi_I$
denotes the submatrix of $\Phi$ obtained by selecting the columns in $I$.

Without loss of generality, $c < 1$.
Since $\Phi_I$ is a $n \times \lceil \d N \rceil$ matrix
and $c_2 n \le \lceil \d N \rceil \le n$, Lemma~\ref{Phi norm} applies
for $\Phi_I$. Taking the union bound over all $I$, we conclude that
for every $t > C_1 \b$
\begin{multline*}
 \P\{ \exists I : \; \|\Phi_I\| > t \sqrt{n} \}
  \le
  \\
  \binom{N}{\lceil \d N \rceil} \exp(-c_1 n t^2/\b^2)
 \le \\
 \exp \big[ (\log(e/\d) - c_1 t^2/\b^2) n \big]
 \le  \\ \exp(-c_3 n t^2/\b^2)
\end{multline*}
if we choose $t = C \b \sqrt{\log \l}$ and use our choice of $\d = c / \l$.
(Here $c_3=c_1/2$ and $C$ are absolute constants).
With this choice of $t$, we can write the estimate above as
\begin{multline*}
\P\{ \exists I : \; \|\Phi_I\| > C \b \sqrt{\frac{\log \l}{\l}} \sqrt{N} \}
\le
\\
\exp(-c_3 C^2 n \log \l)
\le \l^{-n}
\end{multline*}
provided we choose the absolute constant $C$ sufficiently big.
This means that the uncertainty principle with parameters
\eqref{subgaussian parameters} fails with probability at most $\l^{-n}$.
\endproof

\medskip

Unlike random orthogonal or partial Fourier matrices considered
in Sections~\ref{s: orthogonal} and \ref{s: Fourier},
subgaussian matrices do not in general have orthonormal rows.
Nevertheless, the rows of subgaussian matrices are almost orthogonal,
and their columns form almost tight frames as we describe below.
So, one can use Theorem~\ref{almost tight} instead of
Theorem~\ref{algorithm valid} to compute Kashin's representations
for such almost tight frames.

The almost orthogonality of subgaussian matrices can be expressed as follows:

\begin{lemma}                   \label{almost orthogonal}
 Let $\Phi$ be a $n \times N$ matrix whose entries are
 independent mean zero subgaussian random variables with parameter $\b$
 and with variance $1$. There exist constants $C = C(\b)$,
 $c = c(\b) > 0$ such that the following holds.
 Assume that
 $$
 N > \frac{C}{\e^2} \log \Big( \frac{2}{\e} \Big) \cdot n
 $$
 for some $\e \in (0,1)$.
 Then
 $$
 \P \{ \|\id - \frac{1}{N} \Phi \Phi^*\| > \e \}
   \le 2 \exp(-c N \e^2).
 $$
\end{lemma}

\begin{remark}
 The dependence in $C(\b)$, $c(\b)$ is polynomial.
 Explicit bounds can be deduced from \cite{MPT short}.
\end{remark}

\medskip

As a straightforward consequence, we obtain:

\begin{corollary} {\em (Subgaussian frames are almost tight)}
 Let $\Phi$ be a subgaussian matrix as in Lemma~\ref{almost orthogonal}.
 Then the columns of the matrix $\frac{1}{\sqrt{N}} \Phi$
 form an $\e$-tight frame $(u_i)_{i=1}^N$ in $\C^n$.
\end{corollary}

\noindent {\bf Proof of Lemma~\ref{almost orthogonal}. }
In this proof, $C_1, C_2, c_1, c_2, \ldots$ will denote positive
absolute constants.
By a duality argument as in \eqref{duality},
$$
\|\id - \frac{1}{N} \Phi \Phi^*\|
= \sup_{x \in S^{n-1}} |\frac{1}{N} \|\Phi^* x\|_2^2 - 1|.
$$
Denote the columns of $\Phi$ by $\phi_i$.
Fix a vector $x \in \C^n$, $\|x\|_2 = 1$.
Since the entries of the vector $\phi_i$ are subgaussian
with parameter $\b$,
the random variable $\< \phi_i, x\> $ is also subgaussian
with parameter $C_1 \b$, where $C_1$ is an absolute constant
(see Fact~2.1 in \cite{LPRT}).
Moreover, this random variable has mean zero and variance $1$.
We can use Bernstein's inequality (see \cite{VW}) to control
the average of the independent mean zero random variables
$|\< \phi_i, x\> |^2 - 1$ as
$
\P \{ |\frac{1}{N} \|\Phi^* x\|_2^2 - 1| > u \}
=
$
\begin{multline*}
\P \Big\{ \Big| \frac{1}{N} \sum_{i=1}^N |\< \phi_i, x\> |^2 - 1 \Big| > u \Big\}
\le
\\
2 \exp( -c_1 N u^2/\b^4)
\end{multline*}
for all $u \le c\b$, where $c_1>0$ is an absolute constant.

Denote $U = \frac{1}{\sqrt{N}} \Phi$.
There exists a $u$-net $\NN$ of the sphere $S^{n-1}$ in the Euclidean
norm, and with cardinality $|\NN| \le (3/u)^n$ (see e.g. \cite{Pi} Lemma~4.16).
Using the probability estimate above, we can take the union bound to
estimate the probability of the event
$$
A := \{ \forall y \in \NN :\; |\|U^*y\|_2^2 - 1| \le u \}
$$
as
$$
\P(A^c) \le (3/u)^n \cdot 2 \exp( -c_1 N u^2/\b^4).
$$
Applying Lemma ~\ref{Phi norm} with $t = C_1 \b$, we see that the event
$
B := \{ \|U^*\| \le C_1 \b \}$ $
\text{satisfies} \ \ \
\P(B^c) \le \exp(-c_2 N).
$
Consider a realization of the random variables for which the event
$A \cap B$ holds.
For every $x \in S^{n-1}$, we can find an element of the net $y \in \NN$
such that $\|x-y\|_2 \le u$, which implies by the triangle inequality that
\begin{multline*}
|\|U^*x\|_2 - 1|
\le
\\
|\|U^*y\|_2 - 1| + |\|U^*x\|_2 - \|U^*y\|_2|
\le
\\
|\|U^*y\|_2^2 - 1| + \|U^*(x-y)\|_2
\le
\\
u + 2 C_1 \b u \le C_2 \b u,
\end{multline*}
where $C_2 = 1+2C_1$.
Now let $u = \e / 3 C_2 \b$. Thus $C_2 \b u = \e/3 \in (0,1)$, and
the estimate above yields $|\|U^*x\|_2^2 - 1| < \e$ for all $x \in S^{n-1}$
once the event $A \cap B$ holds. Thus
\begin{multline*}
\P \{ \|\id - \frac{1}{N} \Phi \Phi^*\| > \e \}
\le
\\
\P \{ \exists x \in S^{n-1} :\; |\|U^*x\|_2^2 - 1| > \e \}
\le
\\
\P((A \cap B)^c) \le \\
(3/u)^n \cdot 2 \exp( -c_1 N u^2/\b^4) + \exp(-c_2 N)
\le
\\
2 \exp(-c N \e^2)
\end{multline*}
by our choice of $u$ and by the assumption on $N$.
\endproof

\section*{Acknowledgment}
This project started in 2003 when the first author visited University of California, Davis.
The second author thanks the participants of Time-Frequency Brown Bag Seminar
at Princeton for their valuable comments, and especially Ingrid~Daubechies for her
encouragement and interest in this work.
Both authors thank the referees for their many suggestions that helped to greatly
improve the presentation.

\bigskip 

\begin{small}
Yurii Lyubarskii earned the Ph.D. in mathematics from the Institute of Low Temperature Physics and Engineering in Kharkov, Ukraine, in 1974, and D.Sci. in Mathematics from the Leningrad Branch of Steklov Mathematical Institute in 1990. 

During 1971-2001 he was a researcher and a leading researcher at Institute of Low Temperature Physics and Engineering. Since 1996 he is a Professor of Mathematics at  Norwegian University of Science and Technology in Trondheim, Norway. His research interests include complex analysis, harmonic analysis, and applications to signal processing.

\bigskip

Roman Vershynin earned the Ph.D. degree in mathematics from the University of Missouri-Columbia in 2000. 

He was a postdoctoral researcher at Weizmann Institute of Science, Israel, in 2000-2001, and a PIMS postdoctoral researcher at University of Alberta, Canada, in 2001-2003. He was an Assistant and later Associate Professor of Mathematics at University of California, Davis in 2003-2008. Since 2008 he is a Professor of Mathematics at University of Michigan. His research interests include geometric functional analysis, convex geometry, probability theory, and numerical algorithms. 
\end{small}

\end{document}